%% file: main.tex
\documentclass[journal]{IEEEtran}
\pdfoutput=1
\usepackage{amsmath,amsfonts,amsthm,amssymb}
\usepackage[ruled]{algorithm2e}
\usepackage{graphicx}
\usepackage{cases}
\usepackage{cite}
\usepackage{color}
\usepackage{multirow}
\usepackage{booktabs}
\usepackage{bm}
\usepackage{hyperref}
\allowdisplaybreaks[3]
\newtheorem{theorem}{Theorem}
\newtheorem{lemma}{Lemma}
\newtheorem{assumption}{Assumption}
\newtheorem{definition}{Definition}
\newtheorem{proposition}{Proposition}
\newtheorem{corollary}{Corollary}

\hypersetup{hidelinks}

\begin{document}

\title{Approximating Dispatchable Regions in Three-Phase Radial Networks with Conditions for Exact SDP Relaxation}

\author{Bohang Fang, Yue Chen, \IEEEmembership{Member, IEEE},  Changhong Zhao, \IEEEmembership{Senior Member, IEEE}
\thanks{This work was supported by Hong Kong Research Grants Council under Grant No. 14212822 and the National Natural Science Foundation of China under Grant No. 52307144. (Corresponding author: Changhong Zhao.)}
\thanks{B. Fang and C. Zhao are with the Department of Information Engineering, the Chinese University of Hong Kong, New Territories, Hong Kong. Emails: \{fb021, chzhao\}@ie.cuhk.edu.hk}
\thanks{Y. Chen is with the Department of Mechanical and Automation Engineering, the Chinese University of Hong Kong, New Territories, Hong Kong. Email: yuechen@mae.cuhk.edu.hk}
}

% \markboth{IEEE Transactions on Power Systems,~Vol.~xx, No.~xx, April~2023}
% {xxx}

\maketitle

\begin{abstract}
The concept of dispatchable region plays a pivotal role in quantifying the capacity of power systems to accommodate renewable generation. In this paper, we extend the previous approximations of the dispatchable regions on direct current (DC), linearized, and nonlinear single-phase alternating current (AC) models to unbalanced three-phase radial (tree) networks and provide improved outer and inner approximations of dispatchable regions. Based on the nonlinear bus injection model (BIM),
we relax the non-convex problem that defines the dispatchable region to a solvable semidefinite program (SDP) and derive its strong dual problem (which is also an SDP).
Utilizing the special mathematical structure of the dual problem, an SDP-based projection algorithm is developed to construct a convex polytopic outer approximation to the SDP-relaxed dispatchable region.
Moreover, we provide sufficient conditions to guarantee the exact SDP relaxation by adding the power loss as a penalty term, thereby providing a theoretical guarantee for determining an inner approximation of the dispatchable region. 
Through numerical simulations, we validate the accuracy of our approximation of the dispatchable region and verify the conditions for exact SDP relaxation.
\end{abstract}
        

\section*{Nomenclature}
\addcontentsline{toc}{section}{Nomenclature}
\subsection*{Constant Parameters}
\begin{IEEEdescription}
[\IEEEusemathlabelsep\IEEEsetlabelwidth{$i, j, x,k$}] %set the maximum width of the symbol
\item [$\mathbf{Y}$] Network admittance matrix.
\item [$\bm{\mathcal{P}}_{j}^\phi\!,\!\bm{\mathcal{Q}}_{j}^\phi$] Coefficient matrices related to the active and reactive power injections at bus $j$ phase $\phi$.
\item[$\underline{p}_j^\phi, \overline{p}_j^\phi$] Lower and upper limits of controllable active power injection at bus $j$ phase $\phi$.
\item[$\underline{q}_j^\phi, \overline{q}_j^\phi$] Lower and upper limits of controllable reactive power injection at bus $j$ phase $\phi$.
\item[$\underline{v}_j^\phi, \overline{v}_j^\phi$] Lower and upper safety limits of voltage at bus $j$ phase $\phi$.
\item [$\bm{v}_{\textnormal{ref}}$] Voltage of the slack (root) bus.
\item [$\textbf{C}$] Cost matrix related to the active power loss.
\end{IEEEdescription}
\subsection*{Variables}
\begin{IEEEdescription}
[\IEEEusemathlabelsep\IEEEsetlabelwidth{$i, j$}] %set the maximum width of the symbol
\item [$u_j^\phi$] Active power generation of the renewable generator at bus $j$ phase $\phi$.
\item [$z_{kj}^\phi$] The $k$-th ($k=1\sim6$) non-negative slack variable on bus $j$ phase $\phi$. 
\item [$\lambda_{kj}^\phi$] The $k$-th ($k=1\sim6$) non-negative dual variable for inequality constraint on bus $j$ phase $\phi$.
\item [$\bm{\alpha}$] Dual variables for equality constraint.
\item [$\mathbf{W}$] Aggregated squared voltage magnitude.
\item [$\mathbf{A}$] Positive semidefinite dual variable in the FP$'(u)$.
\item [$\mathbf{B}$] Positive semidefinite dual variable in the DE$'(u)$.
\end{IEEEdescription}
\subsection*{Optimization Problems, Values, Sets}
\begin{IEEEdescription}
[\IEEEusemathlabelsep\IEEEsetlabelwidth{$i, j, x, y$}] %set the maximum width of the symbol
\item [FP$(\mathbf{u})$] Feasibility problem for renewable power $\mathbf{u}$.
\item [fp$(\mathbf{u})$] Minimum objective value of FP$(\mathbf{u})$.
\item [$\mathcal{U}$] Dispatchable region of $\mathbf{u}$, in which fp$(\mathbf{u})=0$.
\item [FP$'(\mathbf{u})$] Semidefinite program (SDP) relaxation of FP$(\mathbf{u})$.
\item [fp$'(\mathbf{u})$] Minimum objective value of FP$'(\mathbf{u})$.
\item [$\mathcal{U}'$] SDP-relaxed dispatchable region of $\mathbf{u}$.
\item [DP$'(\mathbf{u})$] Dual problem of FP$'(\mathbf{u})$, also an SDP.
\item [dp$'(\mathbf{u})$] Maximum objective value of DP$'(\mathbf{u})$.
\item [$\mathcal{U}'_{poly}$] Polytopic approximation of $\mathcal{U}'$, returned by Algorithm \ref{Algo:cutting}.
\item [$\mathcal{U}'_{poly}(c)$] Polytopic approximation returned by the $c$-th iteration of Algorithm \ref{Algo:cutting}.
\item[PE$(\mathbf{u})$] Adding power loss as a penalty term to FP$'(\mathbf{u})$.
\item[DE$(\mathbf{u})$] Dual problem of PE$(\mathbf{u})$.
\item [$\tilde{\mathcal{U}}$] Inner approximation of $\mathcal{U}$.
\end{IEEEdescription}

\section{Introduction}
\IEEEPARstart{W}{ith} the high penetration of renewable energy sources (RESs), their volatile and intermittent nature presents significant challenges to the operation of electrical grids, particularly the distribution system \cite{Zeng2024AR_Online}.
The distribution system is characterized by limited controllable units and a stronger coupling of active and reactive power flows, primarily due to its high resistance to reactance ratio \cite{Jia2023highratio}, which makes it more difficult to accommodate the fluctuating renewable energy generation.
Therefore, it is crucial to accurately assess the renewable power capacities that can be safely integrated into a distribution network before its actual operation \cite{Chen2023SOCP}. This requires determining all renewable power outputs that ensure the compliance with \emph{safety limits} and \emph{solvability} of power flow equations.

To address the first requirement of adhering to \emph{safety limits}, two widely studied concepts, namely the Do-not-exceed limit (DNEL) \cite{Zhao2015DNE} and the dispatchable region \cite{Wei2015dispatchable} have been proposed in optimization-based methods.
The DNEL delineates an acceptable power range for each renewable generator through robust optimization. To better enhance the accuracy of DNEL, a data-driven approach \cite{Yanqi2022datadriven}, \cite{Awadalla2024datadriven} has been integrated.
However, this approach overlooks the correlation between different renewable generators and hence potentially obtains conservative capacity regions \cite{Can2018MUB, Zhigang2023Linear, Wei2015dispatchable, Shen2022Region, Nguyen2019innerapprox, Chen2021energy-sharing}.
In contrast, the dispatchable region considers such correlation, providing the ranges of RES power injections that can be accommodated by the power system without operational violations \cite{Zeng2024AR_Online}.
This aspect is critical in the hourly assessment of grid adequacy and security \cite{Hegazy2003adequacy}.
Furthermore, the real-time dispatchable region, as defined in \cite{Zhigang2023Linear}, represents the maximum ranges of power injections that can be handled by the power grid within a specific dispatch interval based on a predetermined base point.

The second requirement is \emph{solvability} of power flow equations.
The Banach fixed-point theorem \cite{Bolognani2016PF, Berstein2018loadflow, Fang2025BFS} and the Brouwer fixed-point theorem \cite{Nguyen2018Brouwer} provide conditions for the existence of power flow solution. However, these methods struggle to handle inequality safety constraints.
Sufficient conditions proposed in \cite{Steven2014exactness, fengyu2019, Bose2015exact} for the exact semidefinite program (SDP) relaxation can guarantee the recovery of power flow solutions that satisfy inequality safety constraints.

Alternating current (AC) power flow is essential for distribution networks. Nonetheless, most existing works, including those mentioned earlier \cite{Zhao2015DNE}, \cite{Wei2015dispatchable}, \cite{Yanqi2022datadriven}, \cite{Zhigang2023Linear}, primarily focus on the direct current (DC) power flow model, which results in an excessively optimistic region containing potentially insecure operating points.
For AC power flow models, linearization methods have been employed to approximate the dispatchable regions \cite{Can2018MUB, Zhigang2023Linear, Chen2021energy-sharing}.
However, the nonlinearity of power system models
becomes pronounced under heavy loading conditions, rendering the linear approximation less accurate.
Based on single-phase AC networks, references \cite{Shen2022Region, Chen2023SOCP} solve second-order cone program (SOCP) relaxations to derive a set of boundary points, each enforcing distinct safety limits, and then heuristically construct a dispatchable region as the convex hull of these boundary points. Nevertheless, their methods are not directly applicable to the more practical three-phase unbalanced networks.

Additionally, the actual dispatchable region is inherently non-convex, not polytopic, yet prior studies often approximate it using polytopes while overlooking its non-convex nature.
Thus, feasibility within the obtained region is not guaranteed.
It is worth noting that sufficient conditions in \cite{Steven2014exactness, fengyu2019, Bose2015exact} for exact SDP relaxation are conservative and may not be applicable in this case.
While the convex-concave procedure (CCP) used in \cite{Shen2022Region} identifies the actual dispatchable region, it might converge solely to a local optimum, lacking theoretical
guarantees as noted in \cite{Lipp2016CCP}.
To ensure feasibility, certified inner approximations of dispatchable regions were derived from convex programs based on a tightened-relaxed SOCP \cite{Nick2018Exact}, or refined linear approximations \cite{Nguyen2019innerapprox}.
However, these approximations typically cater to a given objective function, and hence can only explore the maximum allowable range in a particular direction of the vector space of renewable outputs.

We develop an improved method to supplement the literature above. Our main contributions are two-fold:
\begin{enumerate}
    \item \emph{Outer approximation of dispatchable regions by the SDP relaxation.}
    Extending the results in \cite{Chen2023SOCP, Shen2022Region} from single-phase to unbalanced three-phase networks, we use the bus injection model (BIM) and SDP relaxation to directly formulate the dispatchable region.
    By exploiting the structure of the dual SDP problem, we introduce a projection algorithm (Algorithm \ref{Algo:cutting}) for constructing a polytopic outer approximation of the SDP-relaxed dispatchable region.
    \item \emph{Inner approximation of dispatchable regions by sufficient conditions to guarantee exact SDP relaxation.} 
    By adding the power loss as a penalty term, we provide a rigorous proof that the globally optimal solution to the SDP relaxed problem is exact (i.e., feasible for the original problem before relaxation) under certain practical conditions, thereby finding an inner approximation of the dispatchable region.
    We complete the proof by contradiction and the properties of the $\mathcal{G}$-invertible matrix proposed in \cite{VANDERHOLST2003SDPmatrix}.
\end{enumerate}

The rest of this paper is organized as follows. Section \ref{sec:dispatable region} defines the dispatchable region by an optimization problem and subsequently relaxes it to an SDP. Section \ref{sec:polytopic approximation} elaborates our proposed method for outer approximating the dispatchable region using a polytope. Section \ref{sec:exactness} introduces conditions for the exact SDP relaxation to find an inner approximation of the dispatchable region. Section \ref{sec:numerical} reports numerical experiments and Section \ref{sec:conclusion} concludes the paper.

\section{Dispatchable region and relaxation}\label{sec:dispatable region}
Accommodating a maximum amount of renewable generation while ensuring system security is a primary operational goal in modern power systems.
Traditionally, surplus renewable power is curtailed when the flexibility of power system is exhausted.
However, the curtailment actions are not immediate and are frequently implemented after the system has already encountered security violations.
To address this challenge, a more proactive approach is necessary. 
The concept of the dispatchable region can offer such an approach by representing the area of feasible renewable generations.
By leveraging the dispatchable region, operators can set boundaries for renewable generators to prevent them from surpassing predefined thresholds \cite{Zhao2015DNE}.

\begin{figure}[ht]
\centering
\includegraphics[width=0.8 \columnwidth]{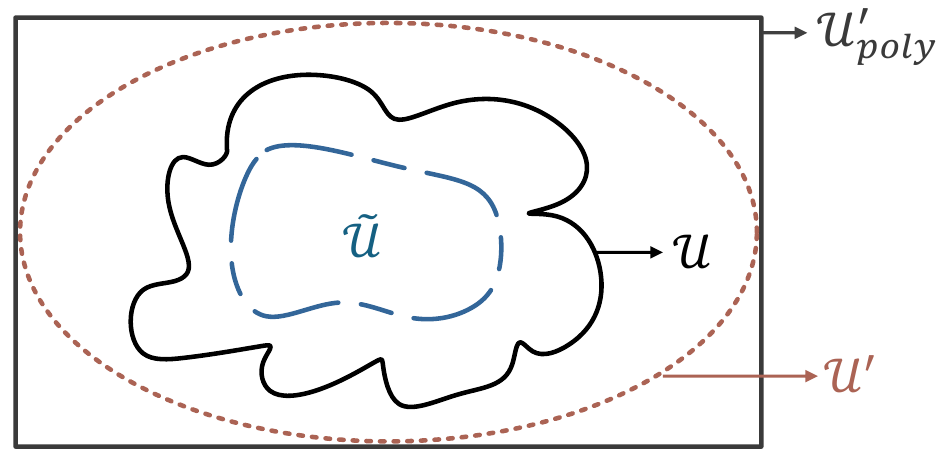}
\caption{Relationship between the dispatchable region and our approximations: $\mathcal{U}$ represents the dispatchable region; $\mathcal{U}'$ is its SDP-relaxed region; $\mathcal{U}'_{poly}$ denotes the outer polytopic approximation; and $\tilde{\mathcal{U}}$ depicts an inner approximation of $\mathcal{U}$.}
\label{fig:DispToy}
\end{figure}

We formulate an optimization problem to define the dispatchable region under the bus injection model (BIM) and then relax it to a convex semidefinite program (SDP). Their relation is shown in Fig. \ref{fig:DispToy}.

\subsection{Power Network Model}
Consider a distribution network modeled as a radial (tree) graph on a set $\mathcal{N}={0,1,...,N}$ of buses (vertices), where $0$ indexed the root bus (slack bus) and $\mathcal{N}^+=\mathcal{N}\backslash 0$. The network has three phases, where each phase of a bus is referred to as a ``node''. 
We use $(j,k)$ and $j\sim k$ interchangeably to denote a line connecting buses $j$ and $k$. 
Each three-phase line $(j,k)$ is characterized by its series admittance matrix $\mathbf{y}_{jk}\in \mathbb{C}^{3\times3}$. 
The network admittance matrix $\mathbf{Y}\in\mathbb{C}^{3(N+1)\times3(N+1)}$ has $(N+1)\times(N+1)$ blocks, each block being a $3\times3$ matrix as follows: 
\begin{align*}
    \mathbf{Y}_{jj}&= \sum_{k: j\sim k} \mathbf{y}_{jk}, \ j \in \mathcal{N} \qquad
    \mathbf{Y}_{jk}&= \begin{cases} -\mathbf{y}_{jk} &, j\sim k\\ 0 &, j\nsim k.
\end{cases}   
\end{align*}
Let the complex voltage at bus $j$, phase $\phi\in\{a,b,c\}$ be $V_j^{\phi}$ and $\mathbf{V}_j = [V_j^a, V_j^b, V_j^c]^{\sf T}$. 
Let $\mathbf{V}=[\mathbf{V}_0^{\sf T}, \mathbf{V}_1^{\sf T}, ..., \mathbf{V}_N^{\sf T}]^{\sf T}$ be the voltage vector of the entire network.
Let $s_j^\phi$ be a controllable power injection into phase $\phi$ at bus $j$, and define $\mathbf{s}_j=[s_j^a, s_j^b, s_j^c]^{\sf T}$ and $\mathbf{s}=[\mathbf{s}_1^{\sf T},...,\mathbf{s}_N^{\sf T}]^{\sf T}$. 
Let $\mathbf{e}_j^{\phi}\in\mathbb{R}^{3(N+1)}$ be the base vector with $1$ at its $\left(3j+\phi\right)$-th entry ($\phi = 1,2,3$ to equivalently represent $a,b,c$) and zero otherwise. 
With $\mathbf{E}_j^{\phi}:=\mathbf{e}_j^{\phi}(\mathbf{e}_j^{\phi})^{\sf T}$, define
$\mathbf{Y}_j^{\phi}\ := \ \mathbf{E}_j^{\phi}\mathbf{Y} \ \in\mathbb{C}^{3(N+1)\times3(N+1)}$.

We further define:
\begin{align*}
    \bm{\mathcal{P}}_j^{\phi}\ := \ \frac{1}{2}\left((\mathbf{Y}_j^\phi)^{\sf H}+\mathbf{Y}_j^{\phi}\right), && 
    \bm{\mathcal{Q}}_j^{\phi}\ := \ \frac{1}{2i}\left((\mathbf{Y}_j^{\phi})^{\sf H}-\mathbf{Y}_j^{\phi}\right),
\end{align*}
where $i$ denotes the imaginary unit and both $\bm{\mathcal{P}}_j^{\phi}$ and $\bm{\mathcal{Q}}_j^{\phi}$ are Hermitian.  
We adopt the three-phase bus injection model (BIM), which compactly expresses the power injections as \cite{fengyu2019}:
\begin{eqnarray}\label{eq:power injection}
    \textnormal{Re}(s_j^\phi)+u_j^\phi \ = \ \mathbf{V}^{\sf H}\bm{\mathcal{P}}_j^{\phi}\mathbf{V},\qquad
    \textnormal{Im}(s_j^{\phi})\ = \ \mathbf{V}^{\sf H}\bm{\mathcal{Q}}_j^{\phi}\mathbf{V}
\end{eqnarray}
where $u_j^\phi$ denotes the renewable active power generation at phase $\phi$ at bus $j$. Suppose a subset $\mathcal{N}_U\subseteq\mathcal{N}^+$ of buses has renewable generation, and $U:=|\mathcal{N}_U|$. All buses $j\notin\mathcal{N}_U$ and all their phases $\phi$ have constant $u_j^\phi\equiv0$.

Assume known capacity limits of controllable power:
\begin{subequations}\label{eq:injection:limit}
    \begin{eqnarray}
         \underline{p}_j^\phi\leq p_j^\phi=\textnormal{Re}(s_j^\phi) \leq \overline{p}_j^\phi,\quad \forall \ j\in\mathcal{N}^+, \  \phi=a,b,c \\
         \underline{q}_j^\phi\leq q_j^\phi=\textnormal{Im}(s_j^\phi) \leq \overline{q}_j^\phi,\quad \forall \ j\in\mathcal{N}^+, \  \phi=a,b,c
    \end{eqnarray}
\end{subequations}
At any node $j$ where there are only fixed (or zero) power injections, the constant limits can be set as $\underline{p}_j^\phi=\overline{p}_j^\phi(=0)$ and/or $\underline{q}_j^\phi=\overline{q}_j^\phi(=0)$. In addition, the following safety limits must be satisfied:
\begin{align}\label{eq:constraint:voltsafety}
    \underline{V}_j^\phi \leq |V_j^\phi| \leq \overline{V}_j^\phi, \quad \forall j\in\mathcal{N}^+, \ \phi\in \{a,b,c\}
\end{align}

\subsection{SDP Relaxation of Dispatchable Region}
With the model above, the formal definition of the dispatchable region is given below.
\begin{definition}\label{def:dispatchable region}(Dispatchable region): A vector of renewable power generation $\mathbf{u}\in\mathbb{R}^{3U}_{\geq0}$ has a feasible dispatch if there exist $\mathbf{s}\in\mathbb{C}^{3N}$ and $\mathbf{V}\in\mathbb{C}^{3(N+1)}$ such that $(\mathbf{u, s, V})$ satisfies power flow equations \eqref{eq:power injection}, capacity limits \eqref{eq:injection:limit}, and safety limits \eqref{eq:constraint:voltsafety}.
The dispatchable region of renewable power generation is defined as:
\begin{align*}
    \mathcal{U}:=\{\mathbf{u}\in\mathbb{R}^{3U}_{\geq0}~|~ \mathbf{u}\ \textnormal{has a feasible dispatch.}\}
\end{align*}
\end{definition}

Consider a three-phase distribution network with only $Y$-connected power sources and loads. We introduce the following optimization problem to check if any given vector $\mathbf{u}$ is in the dispatchable region.
\begin{align*}
\textnormal{min} \qquad &  \sum_{j, \phi} \left(z_{1j}^\phi+z_{2j}^\phi+z_{3j}^\phi+z_{4j}^\phi+z_{5j}^\phi+z_{6j}^\phi \right) \\
\textnormal{over} \qquad & \mathbf{V}\in\mathbb{C}^{3(N+1)},\  z_{kj}^{\phi}\geq 0 \ (k=1,...,6), \ p_j^\phi,\ q_j^\phi
\\\textnormal{s.t.} \qquad &  \forall j\in \mathcal{N}^+,~\phi\in \{a,b,c\}    \nonumber \\
& \textnormal{\eqref{eq:power injection}}, \qquad
\mathbf{V}_0\ =\ \mathbf{V}_{\textnormal{ref}}, \\ 
& p_j^\phi-\overline{p}_j^\phi\leq z_{1j}^{\phi},\qquad -p_j^\phi+\underline{p}_j^\phi\leq z_{2j}^{\phi}, \\ 
& q_j^\phi-\overline{q}_{j}^\phi\leq z_{3j}^{\phi}, \qquad -q_j^\phi+\underline{q}_{j}^\phi\leq z_{4j}^{\phi}, \\
& |V_j^\phi|-\overline{V}_j^\phi\leq z_{5j}^{\phi}, \qquad -|V_j^\phi|+\underline{V}_j^\phi\leq z_{6j}^{\phi},
\end{align*}
where $\mathbf{V}_{\textnormal{ref}}\in\mathbb{C}^3$ denotes the constant voltage at the slack bus, for example, $\mathbf{V}_{\textnormal{ref}}=[1,~e^{\mathrm{j}\frac{2}{3}\pi},~e^{-\mathrm{j}\frac{2}{3}\pi}]^{\sf T}$ in per unit. 
Any slack variable $z_{kj}^\phi$ can be as big as needed to satisfy the corresponding inequality constraint, but only $z_{kj}^\phi=0, \forall k,j,\phi$ can guarantee feasibility in terms of \eqref{eq:power injection}-\eqref{eq:constraint:voltsafety}.

Using $\mathbf{W}:= \mathbf{V}\mathbf{V}^{\sf H}$ as the decision variable, we obtain the following equivalent formulation, where variables listed to the right of the colon denote the corresponding dual variable of each constraint.
\begin{subequations}\label{eqPrimal}
\begin{align}
\textnormal{FP($\mathbf{u}$): min}  & \  \sum_{j, \phi} \left(z_{1j}^\phi+z_{2j}^\phi+z_{3j}^\phi+z_{4j}^\phi+z_{5j}^\phi+z_{6j}^\phi \right) \\
\textnormal{over} \qquad & \mathbf{W}\succeq 0, \quad z_{kj}^{\phi}\geq 0 \ (k=1,...,6) 
\end{align}
\begin{align}
\textnormal{s.t.} \qquad & \forall j\in \mathcal{N}^+,~ \phi\in\{a,b,c\}   \nonumber &  \\ 
& \textnormal{tr}(\bm{\mathcal{P}}_j^\phi \mathbf{W})-u_j^{\phi}-\overline{p}_j^{\phi}\leq z_{1j}^\phi   & : \lambda_{1j}^\phi \label{eqSDP: injection}\\
& -\textnormal{tr}(\bm{\mathcal{P}}_j^\phi \mathbf{W})+u_j^{\phi}+\underline{p}_j^\phi \leq z_{2j}^{\phi} & : \lambda_{2j}^\phi \\
& \textnormal{tr}(\bm{\mathcal{Q}}_j^\phi \mathbf{W})-\overline{q}_j^\phi \leq z_{3j}^\phi  &: \lambda_{3j}^\phi \\
& -\textnormal{tr}(\bm{\mathcal{Q}}_j^\phi \mathbf{W}) +\underline{q}_j^\phi \leq z_{4j}^\phi  & : \lambda_{4j}^\phi \\
& \textnormal{tr}(\mathbf{E}_j^\phi \mathbf{W})-\overline{v}_j^\phi \leq z_{5j}^\phi 
 & :\lambda_{5j}^\phi \\
& -\textnormal{tr}(\mathbf{E}_j^\phi \mathbf{W})+\underline{v}_j^\phi \leq z_{6j}^\phi &: \lambda_{6j}^\phi \label{eqSDP:volt:last}
\\
& [\mathbf{W}]_{00}=\bm{v}_{\textnormal{ref}} \quad & :\bm{\alpha} \label{eqSDP:slack bus}\\
& \textnormal{rank}(\mathbf{W})=1 \label{eqSDP:rank}
\end{align}
\end{subequations}

Here, $\overline{v}_j^\phi=|\overline{V}_j^\phi|^2$, $\underline{v}_j^\phi=|\underline{V}_j^\phi|^2$, $\bm{v}_{\textnormal{ref}}=\mathbf{V}_{\textnormal{ref}} \mathbf{V}_{\textnormal{ref}}^{\sf H}$, and $[\mathbf{W}]_{00}$ stands for
the $3\times3$ matrix block in the top left corner of $\mathbf{W}$. Therefore, denoting the minimum objective value of FP($\mathbf{u}$) as fp($\mathbf{u}$), the dispatchable region in Definition \ref{def:dispatchable region} is equivalently:
\begin{align*}
    \mathcal{U}=\{\mathbf{u}\in\mathbb{R}^{3U}_{\geq0}~|~ \textnormal{fp}(\mathbf{u})=0\}.
\end{align*}
Due to the rank constraint \eqref{eqSDP:rank}, problem FP($\mathbf{u}$) is non-convex and therefore hard to analyze. By removing \eqref{eqSDP:rank}, we relax FP($\mathbf{u}$) to a convex SDP: 
\begin{subequations} \label{eqPrimalSDP}
\begin{align}
\textnormal{FP$'$($\mathbf{u}$):} \qquad & \textnormal{min}\ \sum_{k=1,...,6,\, j,\phi} z_{kj}^\phi \\
\textnormal{over} \qquad & \mathbf{W}\succeq 0, \quad z_{kj}^{\phi}\geq 0 \ (k=1,...,6) \label{eqPrimalSDP:psd}
\\\textnormal{s.t.} \qquad & \forall j\in \mathcal{N}^+,~ \phi\in\{a,b,c\}    \nonumber \\ 
& \eqref{eqSDP: injection}-\eqref{eqSDP:slack bus} \label{eqSDR:constraints}
\end{align}
\end{subequations}
Problem FP$'(\mathbf{u})$ defines an SDP\emph{-relaxed} dispatchable region:
\begin{align*}
    \mathcal{U}':=\{\mathbf{u}\in\mathbb{R}^{3U}_{\geq0}~|~ \textnormal{fp}'(\mathbf{u})=0\}
\end{align*}
where fp$'(\mathbf{u})$ is the minimum objective value of FP$'(\mathbf{u})$. It is obvious that $\mathcal{U}\subseteq \mathcal{U}'$ as shown in Fig. \ref{fig:DispToy}.

One method to simplify the dispatchable region is relaxing the semidefinite constraint $\mathbf{W}\succeq0$ to a polytopic cone by constructing a sufficient number of tangent planes \cite{wang2021polyhedral}.
This transforms FP$'(\mathbf{u})$ to a linear program, based on which the algorithms in \cite{Wei2015dispatchable} and \cite{Chen2021energy-sharing} can derive a convex polytopic outer approximation of $\mathcal{U}'$.
Nonetheless, the accuracy of such approximation in satisfying AC power flow equations may be compromised.

In this paper, we directly tackle FP$'(\mathbf{u})$ through its dual problem to preserve the inherent nonlinearity of the AC power flow model, thereby improving the accuracy of our approximation.

\section{Polytopic approximation algorithm}\label{sec:polytopic approximation}
Although FP$'(\mathbf{u})$ becomes convex after the SDP relaxation, determining the set $\mathcal{U}'$ remains challenging due to the interdependencies between the constraints of FP$'(\mathbf{u})$. 
By exploiting the dual problem of FP$'(\mathbf{u})$, we clearly find the influence of $\mathbf{u}$ on the (dual) objective of FP$'(\mathbf{u})$, which allows us to design a projection method for a polytopic outer approximation of $\mathcal{U}'$.

\subsection{Dual SDP}
Let $(\lambda_{1j}^\phi,...,\lambda_{6j}^\phi)$ denote the dual variables of inequality constraints \eqref{eqSDP: injection}-\eqref{eqSDP:volt:last} in FP$'(\mathbf{u})$, and $\bm{\alpha}$ for \eqref{eqSDP:slack bus}. Then the dual of FP$'(\mathbf{u})$ is also an SDP:
\begin{align}
\textnormal{DP$'(\mathbf{u})$: \ max}  & \sum_{j,\phi}\Big((\lambda_{2j}^\phi-\lambda_{1j}^\phi)u_j^\phi-\lambda_{1j}^\phi \overline{p}_j^\phi+\lambda_{2j}^\phi\underline{p}_j^\phi-\lambda_{3j}^\phi\overline{q}_j^\phi \nonumber \\ &+\lambda_{4j}^\phi\underline{q}_j^\phi-\lambda_{5j}^\phi\overline{v}_j^\phi+\lambda_{6j}^\phi\underline{v}_j^\phi\Big)-\textnormal{tr}(\bm{\alpha} \bm{v}_{\textnormal{ref}})  \nonumber \\
\textnormal{over} \qquad & \lambda_{kj}^{\phi}\geq 0 \ (k=1,...,6), \quad \bm{\alpha} \nonumber
\\ \textnormal{s.t.}  
\qquad & \lambda_{kj}^{\phi}\leq 1, \quad \mathbf{A}(\bm{\lambda, \alpha}) \succeq0  \label{eq:dual:constraint}
\end{align}
where 
\begin{align*}
    \mathbf{A}(\bm{\lambda, \alpha}):=& \bm{\Pi}(\bm{\alpha})+\sum_{j,\phi}\Big((\lambda_{1j}^\phi-\lambda_{2j}^\phi)\bm{\mathcal{P}}_j^\phi+(\lambda_{3j}^\phi-\lambda_{4j}^\phi)\bm{\mathcal{Q}}_j^\phi \nonumber \\ 
    & +(\lambda_{5j}^\phi-\lambda_{6j}^\phi)\mathbf{E}_j^\phi\Big)
\end{align*}
and $\bm{\Pi}(\bm{\alpha})$ is a $3(N+1)\times3(N+1)$ matrix whose top left block is $\bm{\alpha}\in\mathbb{C}^{3\times3}$ and other elements are 0.
The Karush–Kuhn–Tucker (KKT) condition for the optimal solution of FP$'(\mathbf{u})$ and DP$'(\mathbf{u})$ is: 
\begin{align}\label{eq:dual:KKT}
    \begin{cases}
        \textnormal{tr}(\mathbf{A}^{\text{opt}}\mathbf{W}^{\text{opt}})=0 \\
        \sum_{k=1,...,6,\, j,\phi} z_{kj}^{\phi,\text{opt}}(1-\lambda_{kj}^{\phi,\text{opt}})=0
    \end{cases}
\end{align}
where $(\cdot)^{\text{opt}}$ denotes the optimal solution.
Let $D_\mathbf{u}(\bm{\lambda}, \bm{\alpha})$ denote the objective function and dp$'(\mathbf{u})$ denote the maximum objective value of DP$'(\mathbf{u})$. The following propositions lay the foundation for approximating the set $\mathcal{U}'$ via the dual SDP DP$'(\mathbf{u})$.

\begin{proposition}\label{prop:strongdual}
    For all $\mathbf{u}\in\mathbb{R}^{3U}$, strong duality holds between FP$'(\mathbf{u})$ and DP$'(\mathbf{u})$, i.e., their optimal values $\textnormal{fp}'(\mathbf{u})=\textnormal{dp}'(\mathbf{u})$.
\end{proposition}
\begin{IEEEproof}
    Consider an arbitrary $\mathbf{u}\in\mathbb{R}^{3U}$. By adjusting the slack variables $z_{kj}^\phi$ to sufficiently large values, one can always satisfy the inequality constraints in FP$'(\mathbf{u})$. Since FP$'(\mathbf{u})$ is a convex problem, it is sufficient to use Slater's condition \cite[Section 5.2.3]{Boyd2004} to complete the proof.  
\end{IEEEproof}

By Proposition \ref{prop:strongdual}, the relaxed region $\mathcal{U}'$ is equivalent to:
\begin{align}\label{eq:def:u'}
   \mathcal{U}'&=\{\mathbf{u}\in\mathbb{R}^{3U}|\textnormal{dp}'(\mathbf{u})=0\}\nonumber \\
   &=\{\mathbf{u}\in\mathbb{R}^{3U}|D_\mathbf{u}(\bm{\lambda}, \bm{\alpha}) \leq 0, \ \forall(\bm{\lambda}, \bm{\alpha})\ \textnormal{satisfying}\ \eqref{eq:dual:constraint}\}
\end{align}
where the second equality holds because $D_\mathbf{u}(\bm{\lambda}, \bm{\alpha})=0$ can always be attained at the dual feasible point $(\bm{\lambda}, \bm{\alpha})=\mathbf{0}$.
\begin{proposition}\label{prop:convex set}
    $\mathcal{U}'$ is a convex set.
\end{proposition}
\begin{IEEEproof}
    Consider arbitrary $\mathbf{u}_1, \mathbf{u}_2\in\mathcal{U}'$ and $t\in[0, 1]$. Denote $\mathbf{u}_t:=t\mathbf{u}_1+(1-t)\mathbf{u}_2$. Then for every $(\bm{\lambda}, \bm{\alpha})$ satisfying \eqref{eq:dual:constraint}, we have
    \begin{align*}
        D_{\mathbf{u}_t}(\bm{\lambda}, \bm{\alpha})&=tD_{\mathbf{u}_1}(\bm{\lambda}, \bm{\alpha})+(1-t)D_{\mathbf{u}_2}(\bm{\lambda}, \bm{\alpha}) \\
        &\leq t \cdot0+(1-t)\cdot0=0\quad \Rightarrow \quad \mathbf{u}_t \in \mathcal{U}'
    \end{align*}
where the first equality is due to linearity of $D_\mathbf{u}(\bm{\lambda}, \bm{\alpha})$ with respect to $\mathbf{u}$ when $(\bm{\lambda}, \bm{\alpha})$ is fixed, and the inequality holds because $\mathbf{u}_1, \mathbf{u}_2\in\mathcal U'$. By the definition of a convex set, $\mathcal{U}'$ is convex.
\end{IEEEproof}
\subsection{Approximating SDP-relaxed Dispatchable Region}
We propose a polytopic approximation algorithm, i.e., Algorithm \ref{Algo:cutting}, to approximate $\mathcal{U}'$ as defined in \eqref{eq:def:u'}.

The algorithm begins with an initial polytope $\mathcal{U}'_{poly}(0)$ that adequately contains $\mathcal{U}'$. Then, it solves the DP$'(\mathbf{u})$ for each vertex $\mathbf{u}$ of
$\mathcal{U}'_{poly}(0)$ and, for each vertex that violates the definition \eqref{eq:def:u'} of $\mathcal{U}'$, add a plane to truncate $\mathcal{U}'_{poly}(0)$.
Other vertices satisfying \eqref{eq:def:u'} are included in $\mathcal{V}_{safe}$ and do not require further validation. 
This process is repeated until a certain iteration $c$, when all vertices of $\mathcal{U}'_{poly}(c)$ are in $\mathcal{V}_{safe}$, or the maximum iteration count $C_{\textnormal{max}}$ is reached.
As long as the initial $\mathcal{U}'_{poly}(0)$ covers $\mathcal{U}'$,  Algorithm \ref{Algo:cutting} will return a polytopic outer approximation to $\mathcal{U}'$.

\begin{algorithm}[ht]
   \caption{Polytopic Approximation to $\mathcal{U}'$}
   \label{Algo:cutting}
{\bf{Initialization:}} $\mathcal{U}'_{poly}(0)=\{\mathbf{u}\in\mathbb{R}^{3U}_{\geq0}~|~\underline{u}\leq u_j^\phi\leq\overline{u}\}$ for sufficiently low $\underline{u}$ and high $\overline{u}$; $\mathcal{V}_{safe}=\varnothing$; $c=0$, $\textnormal{dp}'_{\textnormal{max}}[0]=99$. A small $\epsilon>0$ as a threshold for being numerically close to zero. A maximum number $C_{\text{max}}$ of iterations.\\
    \While{$\textnormal{dp}'_{\textnormal{max}}[c]>\epsilon$ and $c<C_{\textnormal{max}}$}{
    Obtain vertex set $vert(\mathcal{U}'_{poly}(c))$;\\
    $\mathcal{U}'_{poly}(c+1)=\mathcal{U}'_{poly}(c)$; \
    dp$'_{\text{max}}[c+1]=-1$;\\    
    \For{$\mathbf{u}\in vert(\mathcal{U}'_{poly}(c))$ and $\mathbf{u}\notin \mathcal{V}_{safe}$}{solve DP$'(\mathbf{u})$ to obtain an optimal solution $(\bm{\lambda}_{\text{max}}, \bm{\alpha}_{\text{max}})$ and maximum objective value dp$'(\mathbf{u})$;\\
\eIf{$\textnormal{dp}'(\mathbf{u})\leq\epsilon$}{$\mathcal{V}_{safe}=\mathcal{V}_{safe}\cup \{\mathbf{u}\}$;}{add to $\mathcal{U}'_{poly}(c+1)$ a cutting plane: \\
$\sum_{j,\phi}(\lambda_{\textnormal{max}, 2j}^\phi-\lambda_{\textnormal{max}, 1j}^\phi)u_j^\phi+\sum_{j, \phi}\big(-\lambda_{\textnormal{max}, 1j}^\phi\overline{p}_j^\phi+\lambda_{\textnormal{max}, 2j}^\phi\underline{p}_j^\phi-\lambda_{\textnormal{max}, 3j}^\phi\overline{q}_j^\phi+\lambda_{\textnormal{max}, 4j}^\phi\underline{q}_j^\phi-\lambda_{\textnormal{max}, 5j}^\phi\overline{v}_j^\phi+\lambda_{\textnormal{max}, 6j}^\phi\underline{v}_j^\phi\big)-\textnormal{tr}(\bm{\alpha}_{\textnormal{max}}\bm{v}_{\textnormal{ref}})\leq0$;\\
\If{$\textnormal{dp}'(\mathbf{u})> \textnormal{dp}'_{\textnormal{max}}[c+1]$}{$\textnormal{dp}'_{\textnormal{max}}[c+1]=\textnormal{dp}'(\mathbf{u})$.}
}
}
$c\leftarrow c+1$; }
Return $\mathcal{U}'_{poly}(c)$. 
\end{algorithm}

\begin{theorem}\label{theo:outerapprox}
    The intermediate output $\mathcal{U}'_{poly}(c)$ of Algorithm \ref{Algo:cutting} in an arbitrary iteration $c$ is an outer approximation of $\mathcal{U}'$. The approximation improves as $c$ increases. 
\end{theorem}
\begin{IEEEproof}
    Note that the initial $\mathcal{U}'_{poly}(0)$ contains $\mathcal{U}'$. We next prove that any cutting plane added in Algorithm \ref{Algo:cutting} would not remove any point in $\mathcal{U}'$. To show that, consider an arbitrary $\mathbf{u}$ removed by a cutting plane whose coefficients are $(\bm{\lambda}_{\textnormal{max}}, \bm{\alpha}_{\textnormal{max}})$. Then there must be $D_\mathbf{u}(\bm{\lambda}_{\textnormal{max}}, \bm{\alpha}_{\textnormal{max}})>0$. Since $(\bm{\lambda}_{\textnormal{max}}, \bm{\alpha}_{\textnormal{max}})$ is dual feasible satisfying \eqref{eq:dual:constraint}, we have $\mathbf{u}\notin\mathcal{U}'$ by \eqref{eq:def:u'}. 
    Besides, as $c$ increases, more cutting planes have truncated the output polytope, resulting in an increasingly accurate approximation of $\mathcal{U}'$. 
\end{IEEEproof}

As demonstrated by Theorem \ref{theo:outerapprox}, the SDP-relaxed dispatchable region $\mathcal{U}'$ can be suitably approximated by a finite set of cutting planes.
In contrast to linear approximation models such as those in \cite{Berstein2018loadflow}, \cite{Gupta2021grid-aware}, and the first-order Taylor approximation in \cite{Jabr2019TaylorApprox}, our approach retains the nonlinearity of AC power flow, thus enhancing accuracy of the approximation.
Moreover, the linear approximation methodology in \cite{Berstein2018loadflow} necessitates the determination of two reference operating points, a task that presents considerable challenges due to its dependency on the dispatchable region.

\section{Exactness of SDP relaxation}\label{sec:exactness}
Examining the feasibility of the dispatchable region presents a challenging non-convex task for three reasons.
Firstly, sufficient conditions in \cite{Steven2014exactness, fengyu2019, Bose2015exact} to guarantee the exactness of SDP relaxation do not apply to FP$'(\mathbf{u})$.
Secondly, conventional software packages such as CVX \cite{cvx} always return a solution with maximum rank, and are thus hard to determine the existence of a rank-1 solution. Thirdly, conditions to guarantee the solvability of power flow equations in \cite{Bolognani2016PF, Berstein2018loadflow, Nguyen2018Brouwer} are not suitable for addressing safety (inequality) constraints.

The convex-concave procedure (CCP)-based rank minimization algorithm in \cite{Shen2022Region} iteratively generates a sequence converging to a rank-1 solution. However, this solution may be locally optimal, and the convergence of CCP relies on the initial point, lacking theoretical guarantee as noted in \cite{Lipp2016CCP}.

In this section, we propose conditions to guarantee the exactness of SDP relaxation.
Inspired by \cite{Steven2014exactness}, we add a power loss term into the objective function of FP$'(\mathbf{u})$ to obtain a new SDP denoted as PE$(\mathbf{u})$.
Then we prove the optimal solution to PE$(\mathbf{u})$ is always exact under a practical condition of network conductance (Assumption \ref{assump:admittance}), by contradiction and the properties of the $\mathcal{G}$-invertible matrix. 
The optimal solution to PE$(\mathbf{u})$ is global, with guaranteed convergence -- a theoretical and computational advancement over CCP in both efficiency and robustness.
Finally, we construct an inner approximate dispatchable region $\mathcal{\tilde{U}}$ and prove $\mathcal{\tilde{U}}\subseteq\mathcal{U}$ as shown in Fig. \ref{fig:DispToy}.

\subsection{Primal and Dual SDPs}
With matrix $\mathbf{C}:=\frac{1}{2}(\mathbf{Y}+\mathbf{Y}^{\sf H})$, the total active power loss of the network is $\textnormal{tr}(\mathbf{CW})$. A new SDP PE$(\mathbf{u})$ is formulated as follows: 
\begin{subequations} \label{eqPrimalSDP:powerloss}
\begin{align}
\textnormal{PE($\mathbf{u}$):} \qquad & \textnormal{min}  \quad \beta\cdot \sum_{k=1,...,6, \, j,\phi} z_{kj}^\phi+\textnormal{tr}(\mathbf{CW}) \label{eqPrimalSDP:powerloss:obj}\\
\textnormal{over} \qquad & \mathbf{W}\succeq 0, \quad z_{kj}^{\phi}\geq 0 \ (k=1,...,6)
\\\textnormal{s.t.} \qquad & \forall j\in \mathcal{N}^+,~ \phi\in\{a,b,c\}    \nonumber \\ 
& \eqref{eqSDP: injection}-\eqref{eqSDP:slack bus}
\end{align}
\end{subequations}
where $\beta$ is a factor to weight the sum of slack variables that govern the strictness of constraints \eqref{eqSDP: injection}-\eqref{eqSDP:volt:last}.
As $\beta\rightarrow 0^+$, constraints \eqref{eqSDP: injection}-\eqref{eqSDP:volt:last} are relaxed to permit arbitrarily large slack variables, whereas $\beta\rightarrow +\infty$ enforces hard constraints such that the optimal $z_{kj}^{\phi,\text{opt}}=0$
if one exists.
The exactness proof of PE$(\mathbf{u})$ relies on the selection of $\beta$. The threshold of $\beta$ is intricately associated with the line admittance and topology of the network.
A practical choice is, e.g., $\beta=0.2$.

The dual problem of PE$(\mathbf{u})$ is:
\begin{subequations}\label{eqDualSDP:powerloss}
\begin{align}
\textnormal{DE$(\mathbf{u})$: \ max}  & \sum_{j,\phi}\Big((\lambda_{2j}^\phi-\lambda_{1j}^\phi)u_j^\phi-\lambda_{1j}^\phi \overline{p}_j^\phi+\lambda_{2j}^\phi\underline{p}_j^\phi-\lambda_{3j}^\phi\overline{q}_j^\phi \nonumber \\ &+\lambda_{4j}^\phi\underline{q}_j^\phi-\lambda_{5j}^\phi\overline{v}_j^\phi+\lambda_{6j}^\phi\underline{v}_j^\phi\Big)-\textnormal{tr}(\bm{\alpha} \bm{v}_{\textnormal{ref}})  \\
\textnormal{over} \qquad & \lambda_{kj}^{\phi}\geq 0 \ (k=1,...,6), \quad \bm{\alpha} 
\\\textnormal{s.t.} \qquad & \lambda_{kj}^\phi\leq\beta, \quad \mathbf{B}(\bm{\lambda, \alpha})\succeq0
\end{align}
\end{subequations}
where 
\begin{align*}
\mathbf{B}(\bm{\lambda, \alpha}) :=& \bm{\Pi}(\bm{\alpha})+\sum_{j,\phi}\Big((\lambda_{1j}^\phi-\lambda_{2j}^\phi)\bm{\mathcal{P}}_j^\phi+(\lambda_{3j}^\phi-\lambda_{4j}^\phi)\bm{\mathcal{Q}}_j^\phi \nonumber \\ 
& +(\lambda_{5j}^\phi-\lambda_{6j}^\phi)\mathbf{E}_j^\phi\Big)+\mathbf{C}.
\end{align*}
For notation simplicity, we also denote $\mathbf{B}(\bm{\lambda, \alpha})$ as $\mathbf{B}$ in this section.
The KKT condition for a primal-dual optimal solution of PE$(\mathbf{u})$ and DE$(\mathbf{u})$ is:
\begin{subequations}
\begin{numcases}{}
    \textnormal{tr}(\mathbf{B}^{\text{opt}}\mathbf{W}^{\text{opt}})=0 \label{eq:KKT:tr}   \\
    \sum_{k=1,...,6, \, j, \phi} z_{kj}^{\phi,\text{opt}}(\beta-\lambda_{kj}^{\phi,\text{opt}})=0.
\end{numcases}
\end{subequations}

Based on the optimal solution $z_{kj}^{\phi,\text{opt}}$ to PE$(\mathbf{u})$, we define a region $\tilde{\mathcal{U}}$ as follows:
\begin{align*}
    \tilde{\mathcal{U}}:=\{\mathbf{u}\in\mathbb{R}^{3U}_{\geq0}~\Big|~ \mathbf{u} ~\textnormal{makes} \sum_{k=1,...,6,\,j,\phi}z_{kj}^{\phi,\text{opt}}=0 ~\textnormal{in (\ref{eqPrimalSDP:powerloss})}\}.
\end{align*} 

\subsection{Preliminaries} \label{sec:exact:subsec:preliminary}
In this subsection, we provide results related to the mathematical structure of the positive semidefinite (PSD) matrix $\mathbf{W}$ and the definition of the $\mathcal{G}$-invertible matrix.
These findings build a foundation for the proof of the exact SDP relaxation.

\begin{lemma}\label{lemma:SchurComplement}
\cite[Theorem 4.3]{gallier2010schur} Consider a Hermitian matrix $\mathbf{E}$ in the form of:
\begin{align*}
\mathbf{E}=\begin{bmatrix}
\mathbf{E}_1 & \mathbf{E}_2 \\
\mathbf{E}_2^{\sf H} & \mathbf{E}_3
\end{bmatrix}.    
\end{align*}
Suppose $\mathbf{E}_1\succeq0$ and the pseudo-inverse of $\mathbf{E}_1$ is $\mathbf{E}_1^+$. We have $\mathbf{E}\succeq0$ \emph{if and only if} the column space of $\mathbf{E}_2$ is contained in the column space of $\mathbf{E}_1$ and $\mathbf{E}_3-\mathbf{E}_2^{\sf H}\mathbf{E}_1^+\mathbf{E}_2\succeq0$.
\end{lemma}
Divide any $\mathbf{W}\succeq0$ in the variable space of \eqref{eqPrimalSDP:powerloss} into four blocks:
\begin{align*}
    \mathbf{W}=\begin{bmatrix}
        \mathbf{W}_{00} & \mathbf{W}_{01} \\
        \mathbf{W}_{01}^{\sf H} & \mathbf{W}_{11}
    \end{bmatrix},
\end{align*}
where $\mathbf{W}_{00}\in\mathbb{C}^{3\times3}$, $\mathbf{W}_{01}\in\mathbb{C}^{3\times3N}$, $\mathbf{W}_{11}\in\mathbb{C}^{3N\times3N}$.
With $\mathbf{W}_{00}$, $\mathbf{W}_{01}$, $\mathbf{W}_{11}$, $\mathbf{W}$ treated as $\mathbf{E}_1$, $\mathbf{E}_2$, $\mathbf{E}_3$, $\mathbf{E}$, respectively, by Lemma \ref{lemma:SchurComplement}, we get the following more specific lemma.
\begin{lemma}\label{lemma:PSDstructure}
    Given any $\mathbf{W}\succeq0$ that satisfies $\mathbf{W}_{00}=\mathbf{V}_{\textnormal{ref}}\mathbf{V}_{\textnormal{ref}}^{\sf H}$, i.e., \eqref{eqSDP:slack bus}, there exists a unique (for that $\mathbf{W}$) column vector $\widehat{\mathbf{V}}\in\mathbb{C}^{3N}$ and a unique (for that $\mathbf{W}$) matrix $\mathbf{Q}_W\succeq0$ such that $\mathbf{W}_{01}=\mathbf{V}_{\textnormal{ref}}\widehat{\mathbf{V}}^{\sf H}$ and $\mathbf{W}_{11}=\widehat{\mathbf{V}}\widehat{\mathbf{V}}^{\sf H}+\mathbf{Q}_W$. Moreover, $\textnormal{rank}(\mathbf{W})=1$ if and only if $\mathbf{Q}_W=\mathbf{0}$.
\end{lemma}
To better understand the structure of the PSD matrix $\mathbf{B}(\bm{\lambda, \alpha})$ in DE($\mathbf{u}$), we refer to the following definition \cite{VANDERHOLST2003SDPmatrix}:
\begin{definition}\label{def:G-invertible}
    A PSD matrix $\mathbf{X}$ is $\mathcal{G}$-invertible for a graph $\mathcal{G}(\mathcal{V, E})$ if the following two conditions hold:
    \begin{enumerate}
        \item $\forall(i,j)\in\mathcal{E}$, $[\mathbf{X}]_{ij}$ is invertible.
        \item $\forall i, j \in \mathcal{V}$ such that $i\ne j$ and $(i, j) \notin\mathcal{E}$, $[\mathbf{X}]_{ij}$ is a zero matrix.
    \end{enumerate}
where $[\mathbf{X}]_{ij}$ denotes the block matrix corresponding to bus index $i$ and $j$.\footnote{In this paper, as the root bus is indexed as bus $0$ and we focus on three-phase networks, $[\mathbf{X}]_{ij}$ denotes a $3\times3$ matrix consisting of rows ($3i+1$)--($3i+3$) and columns ($3j+1$)--($3j+3$) of $\mathbf{X}\in\mathbb{C}^{3(N+1)\times3(N+1)}$.}
\end{definition}

We make the following assumption to further guarantee that any feasible $\mathbf{B}(\bm{\lambda, \alpha})$ in DE$(\mathbf{u})$ is $\mathcal{G}$-invertible.
\begin{assumption}\label{assump:admittance}
    For each line $(i, j)\in\mathcal{E}$, the series admittance matrix is symmetric $(\mathbf{y}_{ij}=\mathbf{y}_{ij}^{\sf T})$, and each diagonal element of conductance matrix $\mathbf{g}_{ij}:=\textnormal{Re}\{\mathbf{y}_{ij}\}$ has a much higher magnitude than the off-diagonal elements in the same row.
\end{assumption}
Assumption \ref{assump:admittance} is practical, given that the off-diagonal elements in conductance matrix stem from interactions between phases and are usually much smaller than the conductance of an individual phase. 
\begin{proposition}\label{prop:G-invertible}
    Suppose Assumption \ref{assump:admittance} holds. For sufficiently small $\beta$, matrix $\mathbf{B}$ in the feasible set of \textnormal{DE}$(\mathbf{u})$ is $\mathcal{G}$-invertible. 
\end{proposition}
\begin{IEEEproof}
    For $i\ne j$, and $(i,j)\notin\mathcal{E}$, $[\mathbf{B}]_{ij}=\mathbf{0}_{3\times3}$.
    
    For $(i, j)\in\mathcal{E}$, the block matrix $[\mathbf{B}]_{ij}$ is 
    \begin{align}\label{eq:B:invertible}
       & [\mathbf{B}]_{ij}=\textnormal{diag}(\begin{bmatrix}
            (\lambda_{2i}^{a}-\lambda_{1i}^{a}-1)-i(\lambda_{3i}^a-\lambda_{4i}^a)\\
            (\lambda_{2i}^{b}-\lambda_{1i}^b-1)-i(\lambda_{3i}^b-\lambda_{4i}^b) \\
            (\lambda_{2i}^c-\lambda_{1i}^c-1)-i(\lambda_{3i}^c-\lambda_{4i}^c)
        \end{bmatrix})\frac{1}{2}\mathbf{y}_{ij}\nonumber \\
        & +\frac{1}{2}\mathbf{y}_{ij}^{\sf H}\ \textnormal{diag}(\begin{bmatrix}
            (\lambda_{2j}^{a}-\lambda_{1j}^{a}-1)+i(\lambda_{3j}^a-\lambda_{4j}^a)\\
            (\lambda_{2j}^{b}-\lambda_{1j}^b-1)+i(\lambda_{3j}^b-\lambda_{4j}^b) \\
            (\lambda_{2j}^c-\lambda_{1j}^c-1)+i(\lambda_{3j}^c-\lambda_{4j}^c)
        \end{bmatrix}).
    \end{align}
Given $\lambda_{kj}^{\phi}\leq\beta<<1 \ (k=1,...,4)$, we have $[\mathbf{B}]_{ij}\approx-\frac{1}{2}\mathbf{y}_{ij}-\frac{1}{2}\mathbf{y}_{ij}^{\sf H}=-\mathbf{g}_{ij}$, which by Assumption \ref{assump:admittance} implies that $[\mathbf{B}]_{ij}$ is strictly diagonally dominant.
By Gershgorin's circle theorem, $[\mathbf{B}]_{ij}$ is invertible.
\end{IEEEproof}

Reference \cite{VANDERHOLST2003SDPmatrix} considered matrices in which the non-zero off-diagonal real entries align with the lines $\mathcal{E}$ in $\mathcal{G}(\mathcal{V, E})$. We extend the matrices in \cite{VANDERHOLST2003SDPmatrix} to complex matrices and prove the following result.

\begin{proposition}\label{prop:connected}
    Consider a non-zero vector $\bm{\nu}\in\mathbb{C}^{3(N+1)}$. Let $\Omega(\bm{\nu})\subseteq\{0,...,N\}$ collect the bus indices corresponding to the non-zero $3\times1$ subvectors of $\bm{\nu}$, and let $|\Omega(\bm{\nu})|$ represent the cardinality of $\Omega(\bm{\nu})$.
    Suppose $\mathbf{X}$ is a $\mathcal{G}$-invertible matrix that satisfies $\mathbf{X}\bm{\nu}=\mathbf{0}$ with the smallest $|\Omega(\bm{\nu})|$. Then each pair of buses in $\Omega(\bm{\nu})$ are connected by a line in $\mathcal{G}(\mathcal{V, E})$.
\end{proposition}
\begin{IEEEproof}
    If not, then assume $\Omega(\bm{\nu})=\Omega_1\cup\Omega_2\cdots\cup\Omega_M$ with $M\geq2$ and $\Omega_m$, $m=1,...,M$ nonempty and mutually exclusive, so that each pair of buses in $\Omega(\bm{\nu})$ are connected by a line in $\mathcal{G}(\mathcal{V, E})$ if and only if they are in the same subset $\Omega_m$, $m=1,...,M$. Construct $\tilde{\bm{\nu}}$ in the following manner:  
    \footnote{Subvector $[\bm{\nu}]_k$ only has the $(3k+1)$--$(3k+3)$-th elements of $\bm{\nu}$ corresponding to bus index $k=0,1,...,N$.}
    \begin{align*}
        [\tilde{\bm{\nu}}]_k=\begin{cases}
            [\bm{\nu}]_k  &,k\in\Omega_1 \\
            \mathbf{0} &, k\notin \Omega_1
        \end{cases}.
    \end{align*}
    Then by the $\mathcal{G}$-invertibility of $\mathbf{X}$, we have for all $j\in\Omega_1$ that:
    \begin{align*}        [\mathbf{X}\tilde{\bm{\nu}}]_j &= [\mathbf{X}]_{jj}[\tilde{\bm{\nu}}]_j+\sum_{k: k\sim j}
    [\mathbf{X}]_{jk}[\tilde{\bm{\nu}}]_k \\ 
    &= [\mathbf{X}]_{jj}[\bm{\nu}]_j+\sum_{k: k\sim j}
    [\mathbf{X}]_{jk}[\bm{\nu}]_k = [\mathbf{X}\bm{\nu}]_j=\mathbf{0}.
    \end{align*}
    because buses $k\sim j$, $k\in\Omega(\bm{\nu})$ satisfy $k\in\Omega_1$ and $[\tilde{\bm{\nu}}]_k=[\bm{\nu}]_k$; and buses $k\sim j$, $k\notin\Omega(\bm{\nu})$ satisfy $[\tilde{\bm{\nu}}]_k=\mathbf{0}=[\bm{\nu}]_k$.
    Therefore, we have:
    \begin{align*}
        \tilde{\bm{\nu}}^{\sf H}\mathbf{X}\tilde{\bm{\nu}} = \sum_{j\in\Omega_1}[\tilde{\bm{\nu}}]_j^{\sf H}[\mathbf{X}\tilde{\bm{\nu}}]_j = \sum_{j\in\Omega_1}[\tilde{\bm{\nu}}]_j^{\sf H}\mathbf{0} = 0.
    \end{align*}
By Definition \ref{def:G-invertible}, $\mathcal{G}$-invertible matrix $\mathbf{X}$ is PSD, and then by Cholesky decomposition, there exists a matrix $\mathbf{X}_1$ such that $\mathbf{X}=\mathbf{X}_1^{\sf H} \mathbf{X}_1$. Consequently,
\begin{align*}
    0=\tilde{\bm{\nu}}^{\sf H}\mathbf{X}_1^{\sf H}\mathbf{X}_1\tilde{\bm{\nu}}=||\mathbf{X}_1\tilde{\bm{\nu}}||_2^2 \ \Leftrightarrow \ \mathbf{X}_1\tilde{\bm{\nu}}=\mathbf{0} \ \Leftrightarrow \ \mathbf{X}\tilde{\bm{\nu}}=\mathbf{0}. 
\end{align*}
    
    As $|\Omega(\tilde{\bm{\nu}})|=|\Omega_1|<|\Omega(\bm{\nu})|$ and $\tilde{\bm{\nu}}$ is non-zero by construction, it contradicts the minimality of $\Omega(\bm{\nu})$. Then there must be $M=1$, i.e., each pair of buses in $\Omega(\bm{\nu})$ are connected by a line in $\mathcal{G}(\mathcal{V, E})$.
\end{IEEEproof}

\subsection{Proof of Exact SDP Relaxation}
\label{sec:exact:subsec:SDP}
Based on the preliminaries in Section \ref{sec:exact:subsec:preliminary}, we prove the following theorem that guarantees the SDP relaxation PE$(\mathbf{u})$ is exact. As a result, we obtain an inner approximation of the dispatchable region.

\begin{theorem}\label{theorem:exact}
Under Assumption \ref{assump:admittance}, the optimal solution to the SDP PE$(\mathbf{u})$ \eqref{eqPrimalSDP:powerloss} satisfies the rank constraint \eqref{eqSDP:rank} and can recover a vector of nodal voltages.
\end{theorem}
\begin{IEEEproof}
If not, then there exists an optimal solution $\mathbf{W}^{\text{opt}}$ to \eqref{eqPrimalSDP:powerloss} with $\textnormal{rank}(\mathbf{W}^{\text{opt}})\geq2$.\footnote{Note that $\textnormal{rank}(\mathbf{W}^{\text{opt}})$ cannot be 0 due to the constraint $[\mathbf{W}^{\text{opt}}]_{00}=\bm{v}_{\textnormal{ref}}$.}
Let $\mathbf{B}^{\text{opt}}$ denote the optimal value of $\mathbf{B}(\bm{\lambda, \alpha})$ in the dual problem DE$(\mathbf{u})$ \eqref{eqDualSDP:powerloss}.

Suppose the eigen-decomposition of $\mathbf{W}^{\text{opt}}$ is
\begin{align*}
    \mathbf{W}^{\text{opt}}=\mathbf{V}_1\mathbf{V}_1^{\sf H}+\sum_{r=2}^{R} \rho_r \mathbf{V}_r \mathbf{V}_r^{\sf H}, \quad \mathbf{V}_1=\begin{bmatrix}
        \mathbf{V}_{\textnormal{ref}} \\
        \widehat{\mathbf{V}}
    \end{bmatrix}\in\mathbb{C}^{3(N+1)}
\end{align*}
where $R\geq2$, $\rho_2\geq...\geq\rho_{R}>0$, and $\mathbf{V}_r$ is the normalized eigenvector corresponding to the eigenvalue $\rho_r$.

In the KKT condition \eqref{eq:KKT:tr}, as $\mathbf{B}^{\text{opt}}$ is PSD, we have
\begin{align*}
    0&=\textnormal{tr}(\mathbf{B}^{\text{opt}}\mathbf{W}^{\text{opt}})=\textnormal{tr}(\mathbf{B}^{\text{opt}}\mathbf{V}_1\mathbf{V}_1^{\sf H})+\textnormal{tr}\left(\mathbf{B}^{\text{opt}}\sum_{r=2}^{R}\rho_r \mathbf{V}_r\mathbf{V}_r^{\sf H}\right) \\   &
    =\textnormal{tr}(\mathbf{V}_1^{\sf H}\mathbf{B}^{\text{opt}}\mathbf{V}_1)+\sum_{r=2}^{R}\rho_r\textnormal{tr}(\mathbf{V}_r^{\sf H}\mathbf{B}^{\text{opt}}\mathbf{V}_r)\geq0.
\end{align*}
The first equality holds only when $\mathbf{B}^{\text{opt}}\mathbf{V}_r=\mathbf{0}$ for all $r\leq R$.
By Lemma \ref{lemma:PSDstructure}, the first three elements of $\mathbf{V}_r$ for $r\geq2$ are zeros. Therefore, $|\Omega(\mathbf{V}_r)|\leq N$ for all $r\geq2$. 

We always have a non-zero vector $\tilde{\mathbf{V}}$ that satisfies $\mathbf{B}^{\text{opt}}\tilde{\mathbf{V}}=\mathbf{0}$ with the smallest $|\Omega(\tilde{\mathbf{V}})|$. It follows that $|\Omega(\tilde{\mathbf{V}})|\leq|\Omega(\mathbf{V}_2)|\leq N$. So there exists at least one bus not in $\Omega(\tilde{\mathbf{V}})$.

Suppose bus $j$ is not in $\Omega(\tilde{\mathbf{V}})$ but is connected by a line to a bus $k\in\Omega(\tilde{\mathbf{V}})$. 
Bus $k$ must be the only bus in $\Omega(\tilde{\mathbf{V}})$ that connects to $j\notin\Omega(\tilde{\mathbf{V}})$, because otherwise Proposition \ref{prop:connected} would induce a circle that contradicts the tree topology of the distribution network. Then
\begin{align*}
    [\mathbf{B}^{\text{opt}} \tilde{\mathbf{V}}]_j&=[\mathbf{B}^{\text{opt}}]_{jj}[\tilde{\mathbf{V}}]_j+\sum_{l:l\sim j}[\mathbf{B}^{\text{opt}}]_{jl}[\tilde{\mathbf{V}}]_l \\
    &=[\mathbf{B}^{\text{opt}}]_{jj}\mathbf{0}+[\mathbf{B}^{\text{opt}}]_{jk}[\tilde{\mathbf{V}}]_k+\sum_{l:l\sim j,l\notin\Omega(\tilde{\mathbf{V}})}\!\!\!\!\!\![\mathbf{B}^{\text{opt}}]_{jl}[\tilde{\mathbf{V}}]_l \\
    & = [\mathbf{B}^{\text{opt}}]_{jk}[\tilde{\mathbf{V}}]_k.
\end{align*}
The last equality results from 
$[\tilde{\mathbf{V}}]_l=\mathbf{0}$ for $l\notin\Omega(\tilde
{\mathbf{V}})$. By Proposition \ref{prop:G-invertible}, for $(j,k)\in\mathcal{E}$, the $\mathcal{G}$-invertibility of $\mathbf{B}^{\text{opt}}$ ensures that $[\mathbf{B}^{\text{opt}}]_{jk}$ is invertible. As bus $k\in\Omega(\tilde{\mathbf{V}})$, then $[\tilde{\mathbf{V}}]_k\neq \mathbf{0}$. Consequently, $[\mathbf{B}^{\text{opt}}\tilde{\mathbf{V}}]_j=[\mathbf{B}^{\text{opt}}]_{jk}[\tilde{\mathbf{V}}]_k$ must be non-zero, which contradicts $\mathbf{B}^{\text{opt}}\tilde{\mathbf{V}}=\mathbf{0}$. Therefore, $\mathbf{W}^{\text{opt}} = \mathbf{V}_1\mathbf{V}_1^{\sf H}$, where $\mathbf{V}_1$ is a vector of nodal voltages.
\end{IEEEproof}
\textit{Remark:} 
The invertibility of block matrix $[\mathbf{B}^{\text{opt}}]_{ij}$, $\forall(i, j)\in\mathcal{E}$ is important in the above proof. In \eqref{eq:B:invertible}, $\lambda_{1i}^{\phi, \text{opt}}-\lambda_{2i}^{\phi, \text{opt}}+1$ and $\lambda_{3i}^{\phi, \text{opt}}-\lambda_{4i}^{\phi, \text{opt}}$ are indeed the distribution locational marginal prices (DLMPs) of active and reactive power injections at bus $i$ phase $\phi$ \cite{Anthony2018DLMP}. Particularly, the constant number 1 originates from the power loss term in \eqref{eqPrimalSDP:powerloss:obj} and dominates the DLMPs if the weighting factor $\beta$ is small enough. This makes the DLMPs of different phases close to each other, thus reducing congestion and facilitating the exact SDP relaxation \cite{Steven2014exactness, fengyu2019}.

Theorem \ref{theorem:exact} implies the following corollary.
\begin{corollary}\label{coro:inner:approx}
    $\tilde{\mathcal{U}}$ is an inner approximation of $\mathcal{U}$.
\end{corollary}
\begin{IEEEproof}
Consider $\mathbf{u}\in\tilde{\mathcal{U}}$, then $\sum_{k=1,...,6,\,j,\phi}z_{kj}^{\phi,\text{opt}}=0$ for PE$(\mathbf{u})$ \eqref{eqPrimalSDP:powerloss}. By Theorem \ref{theorem:exact}, the optimal solution to \eqref{eqPrimalSDP:powerloss} is feasible for FP$(\mathbf{u})$. Then we have fp$(\mathbf{u})=0$ and thus $\mathbf{u}\in\mathcal{U}$.
\end{IEEEproof}
Theorem \ref{theo:outerapprox} in Section \ref{sec:polytopic approximation} and Corollary \ref{coro:inner:approx} here verify $\tilde{\mathcal{U}}\subseteq\mathcal{U}\subseteq\mathcal{U}'\subseteq\mathcal{U}'_{poly}$ as shown in Fig. \ref{fig:DispToy}.

\section{Numerical Experiments}\label{sec:numerical}
In this section, we conduct numerical experiments on the IEEE 123-bus network. The impact of several factors is tested, and the proposed method is compared with other approaches to show its advantages.
The programs are run on a desktop computer with the 11th Generation Intel(R) Core(TM) i7-11700 at 2.50 GHz processors and 16GB RAM.
\subsection{Baseline}\label{sec:numerical:subsec:baseline}
Consider three renewable generators located at nodes 23a, 67b, and 35c, labeled as $(u_1, u_2, u_3)$, respectively and six controllable generators, whose locations and lower and upper bounds of active and reactive power injections $(\underline{p}_j^\phi, \overline{p}_j^\phi, \underline{q}_j^\phi, \overline{q}_j^\phi)$ are listed in Table \ref{tab:para:gen} (the Baseline case).
Power injections at all remaining non-slack buses remain fixed, with nodal voltage magnitude limits set to [0.9, 1.1] p.u.

\begin{figure}[ht]
\centering
\includegraphics[width=0.8\columnwidth]{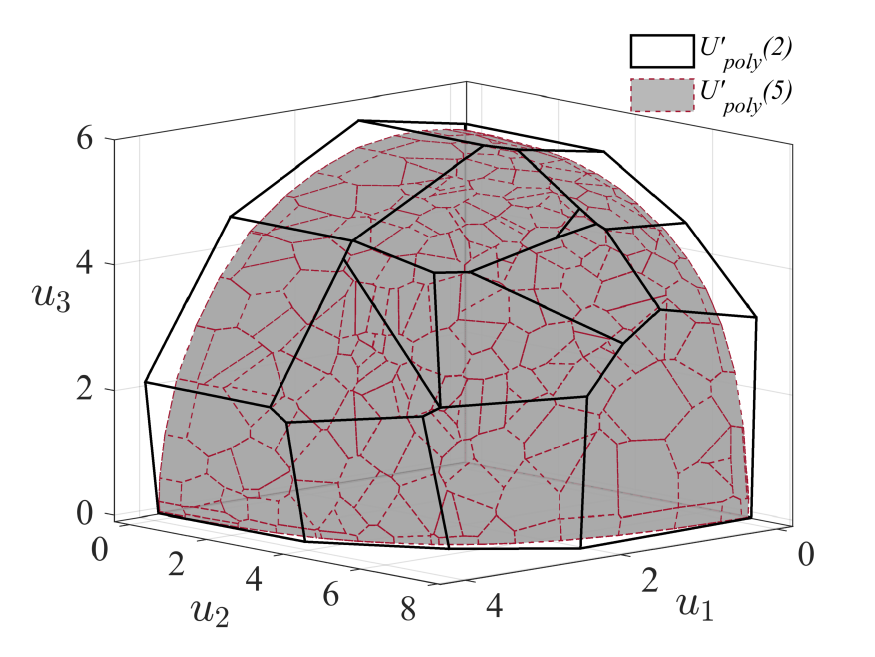}
\caption{Visualization of the convergence behavior of the SDP-relaxed region $\mathcal{U}'$: The gray polytope $\mathcal{U}'_{poly}(5)\approx\mathcal{U}'$ is enclosed by $\mathcal{U}'_{poly}(2)$.}
\label{fig:DP:cmp}
\end{figure}

\begin{figure}[ht]
\centering
\includegraphics[width=0.95 \columnwidth]{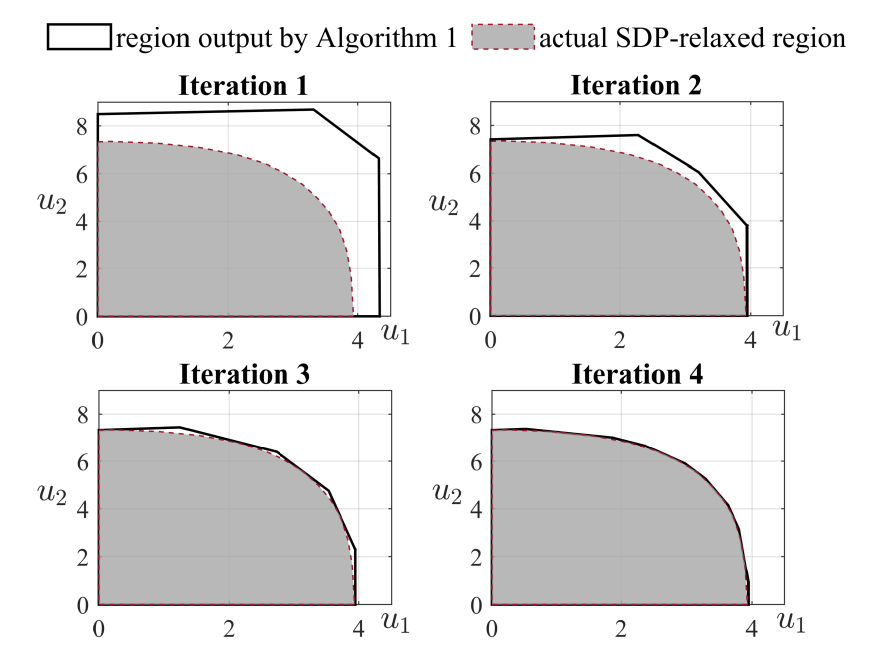}
\caption{Output of Algorithm \ref{Algo:cutting} at different iterations. The polytope $\mathcal{U}'_{poly}(c)$ (solid lines) progressively approximates the SDP-relaxed region $\mathcal{U}'$ (dashed lines).}
\label{fig:poly:iter}
\end{figure}

\begin{table*}[ht]
    \caption{Generator parameters for Baseline, Case L, and Case H scenarios (in p.u.)}
    \centering
    \renewcommand\arraystretch{1.2}
    \begin{tabular}{cl|cccccc|cccccc|cccccc}
    \toprule
    Gen. & Location & \multicolumn{6}{c|}{Baseline}&\multicolumn{6}{c|}{Case L} & \multicolumn{6}{c}{Case H}\\
    No. & &$\underline{p}_i^a$  & $\underline{p}_i^b$  &$\underline{p}_i^c$ &$\overline{p}_i^a$ & $\overline{p}_i^b$& $\overline{p}_i^c$ & $\underline{p}_i^a$  & $\underline{p}_i^b$  &$\underline{p}_i^c$ &$\overline{p}_i^a$ & $\overline{p}_i^b$& $\overline{p}_i^c$ &$\underline{p}_i^a$  & $\underline{p}_i^b$  &$\underline{p}_i^c$ &$\overline{p}_i^a$ & $\overline{p}_i^b$& $\overline{p}_i^c$  \\ 
    \hline
    G1 & bus 13& 0.1& 0.1& 0.1& 0.3& 0.3& 0.3 & 0.1& 0.1& 0.1& 0.2& 0.18& 0.23& 0& 0& 0& 0.3& 0.35& 0.3\\
   G2 & bus 51& 0.2& 0.2& 0.2& 0.5& 0.5& 0.5 & 0.2& 0.2& 0.2& 0.3& 0.35& 0.4& 0& 0& 0& 0.55& 0.5& 0.6\\
   G3 & bus 57& 0.1& 0.1& 0.1& 0.5& 0.5& 0.5& 0.1&0.1 & 0.1& 0.25& 0.3& 0.25& 0& 0& 0& 0.5& 0.6& 0.55\\
   G4 & bus 78& 0.15& 0.15& 0.15& 0.5& 0.5& 0.5&0.15 &0.15 & 0.15& 0.3& 0.2& 0.3& 0& 0& 0& 0.5& 0.6& 0.5\\
   G5 & bus 87& 0.1& 0.1& 0.1& 0.4& 0.4& 0.4& 0.1&0.1 & 0.1& 0.25& 0.25& 0.25& 0& 0& 0& 0.5& 0.45& 0.4\\
   G6 & bus 101& 0.15& 0.15& 0.15& 0.4& 0.4&0.4& 0.15& 0.15& 0.15& 0.25& 0.3& 0.3& 0& 0& 0& 0.4& 0.4&0.4 \\
    \bottomrule
    \end{tabular}   
    \label{tab:para:gen}
\end{table*}
Algorithm \ref{Algo:cutting} is applied with a convergence tolerance of $\epsilon=0.02$ and up to $C_{\textnormal{max}}=6$ iterations. Performance results are summarized in Tables \ref{tab:A1stats} and \ref{tab:A1comp:time} and Figs. \ref{fig:DP:cmp} and \ref{fig:poly:iter}.
Specifically, Algorithm \ref{Algo:cutting} converges after 5 iterations with $\textnormal{dp}'_{\textnormal{max}}[5]=0.0166$.
A detailed statistics of the returned polytopes in each iteration, demonstrating the convergence trend of Algorithm \ref{Algo:cutting} is shown in Table \ref{tab:A1stats}.
The data show a consistent decrease in both the maximum and mean objective values at the polytope vertices over iterations, validating Theorem \ref{theo:outerapprox}. However, there is a trade-off between computational efficiency and polytopic approximation accuracy, as the rates of polytope volume reduction and decrease in $\textnormal{dp}'_{\textnormal{max}}$ slow down while the computational burden increases due to a growing number of vertices.

\begin{table}[ht]   
    \centering
    \renewcommand\arraystretch{1.4}
\caption{Polytope statistics across iterations computed by Algorithm \ref{Algo:cutting}}
    \begin{tabular}{cccccc}
    \toprule
     Iter. & Vertices &  Max & Mean & Volume & Faces\\
    \hline
     0& 8 & 19.14& 9.07& 2208&6 \\
    1&  16& 2.11& 1.27& 229.8& 10\\
    2&  32& 0.443& 0.290& 125.2& 18\\
    3&  62& 0.195& 0.103& 105.1& 33\\
    4&  180& 0.0541& 0.0252& 99.98& 92\\ 
    5&  470& 0.0166& 0.00676& 98.55& 237\\
    \bottomrule
    \end{tabular}
    \label{tab:A1stats}
\end{table}

The tightening of the polyhedral approximation from the second to the fifth iteration is illustrated in Fig. \ref{fig:DP:cmp}.
To enhance visualization, a cross section of the $u_1$-$u_2$ plane (with $u_3=0$) is shown in Fig. \ref{fig:poly:iter}, illustrating that the polytope $\mathcal{U}'_{poly}(4)$ closely approximates the SDP-relaxed region $\mathcal{U}'$ by the fourth iteration.
As the approximation achieves sufficient accuracy at the fourth iteration, further computation would bring negligible improvement. Therefore, the convergence tolerance $\epsilon$ can be increased to improve computational efficiency while maintaining high accuracy.

We set $u_3=0$ and run Algorithm \ref{Algo:cutting} with two renewable generators $(u_1, u_2)$.
The computation time for each iteration is reported in Table \ref{tab:A1comp:time}. Over the first four iterations, it increases as the number of vertices on the polytope grows. 
However, in the fifth iteration, the time decreases, as 9 out of 11 vertices already belong to the set $\mathcal{V}_{safe}$, reducing the need for re-computation.

\begin{table}[ht]   
    \centering
    \renewcommand\arraystretch{1.4}
\caption{Computation time per iteration of Algorithm \ref{Algo:cutting} with two renewable generators}
    \begin{tabular}{cccccc}
    \toprule
     Iter. & 1 &  2 & 3 & 4 & 5\\
    \hline
     Time (s) & 212.3 & 280.8& 284.7& 286.7&178.7 \\
    \bottomrule
    \end{tabular}
    \label{tab:A1comp:time}
\end{table}

\subsection{Impact of Different Factors}\label{sec:numerical:subsec:factors}

\begin{figure}[ht]
\centering
\includegraphics[width=0.95 \columnwidth]{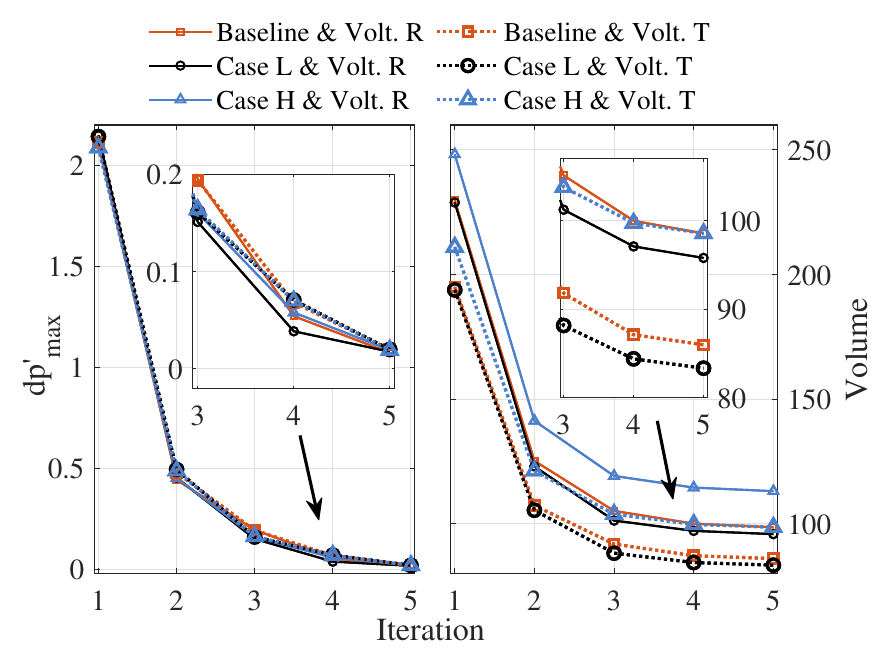}
\caption{Evolution of $\textnormal{dp}'_{\textnormal{max}}$ and polytope volume by Algorithm \ref{Algo:cutting} over iterations for six test cases. Solid lines represent relaxed voltage limits (Volt. R), and dotted lines represent tightened voltage limits (Volt. T). Colors indicate different generator capacities.}
\label{fig:dpmax:volume}
\end{figure}

\begin{figure}[ht]
\centering
\includegraphics[width=0.85 \columnwidth]{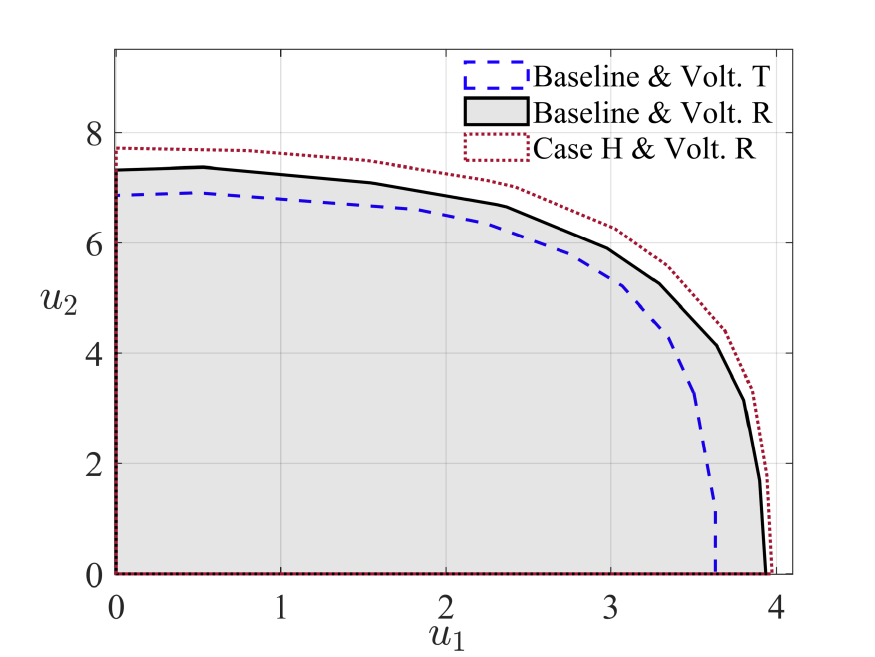}
\caption{Comparison of $\mathcal{U}'_{poly}$ under different scenarios. The solid black line indicates Baseline $\&$ Volt. R, which serves as the reference. The blue dashed and red dotted lines represent Baseline $\&$ Volt. T and Case H $\&$ Volt. R, respectively.}
\label{fig:poly:comp:2d}
\end{figure}

Three cases of generator capacities $[\underline{p}_j^\phi, \overline{p}_j^\phi]$ are considered: (1) \emph{Baseline}, which follows the same setting as in Section \ref{sec:numerical:subsec:baseline}; (2) \emph{Case L}, with reduced capacities, (3) \emph{Case H}, with increased capacities. The parameter details for these cases are listed in Table \ref{tab:para:gen}. Additionally, two voltage safety limits $[ \underline{V}_j^\phi, \overline{V}_j^\phi]$ are examined: (1) \emph{Volt. R}, with relaxed limits of [0.9, 1.1] p.u., and (2) \emph{Volt. T}, with tightened limits of [0.95, 1.05] p.u. 
Combining these parameters yields six test cases, whose results are summarized in Figs. \ref{fig:dpmax:volume} and \ref{fig:poly:comp:2d}, and Table \ref{tab:A1time:3cases}.

The evolution of the maximum vertex objective values $\textnormal{dp}'_{\textnormal{max}}$ and polytope volumes across iterations are illustrated in Fig. \ref{fig:dpmax:volume}. While the values of $\textnormal{dp}'_{\textnormal{max}}$ are similar across different cases, noticeable variations are observed in the volumes of polytopes. Reduced generator capacities or more restrictive voltage limits result in a smaller SDP-relaxed dispatchable region.

To further illustrate these variations, three representative cases: \emph{Baseline $\&$ Volt. T}, \emph{Baseline $\&$ Volt. R}, and \emph{Case H $\&$ Volt. R} are compared in Fig. \ref{fig:poly:comp:2d}. The Baseline $\&$ Volt. R case is selected as the benchmark since the other two cases differ from it in specific ways—one has a tightened voltage limit (resulting in a smaller polytope), while the other has a increased generator capacity (leading to a larger polytope). Their differences are depicted in Fig. \ref{fig:poly:comp:2d}.

Table \ref{tab:A1time:3cases} compares the computation time for the selected cases. While a larger dispatchable region generally requires more computational time, the differences remain modest.

\begin{table}[ht]   
    \centering
    \renewcommand\arraystretch{1.4}
\caption{Computation time for selected cases}
    \begin{tabular}{cccc}
    \toprule
     Case & Baseline \!$\&$\!\! Volt. T &  Baseline \!$\&$\!\! Volt. R & Case H \!$\&$\!\! Volt. R \\
    \hline
     Time (s) & 1094.1 &  1186.9 &  1188.3\\
    \bottomrule
    \end{tabular}
    \label{tab:A1time:3cases}
\end{table}

\subsection{Comparison With Other Methods}
\begin{figure}[ht]
\centering
\includegraphics[width=0.85 \columnwidth]{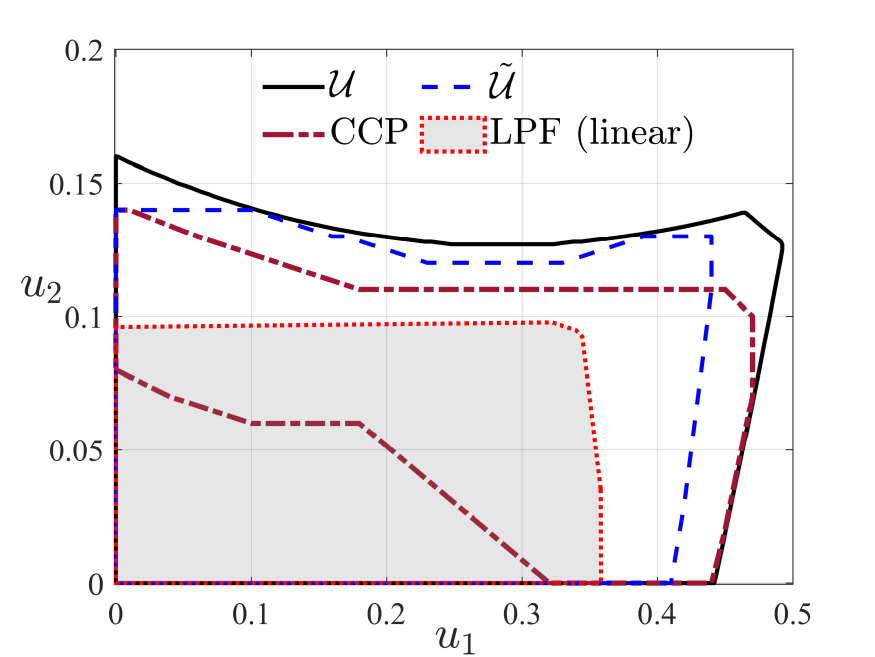}
\caption{Comparison of dispatchable region approximations: solid black line delineates the dispatchable region $\mathcal{U}$; blue dashed line indicates $\tilde{\mathcal{U}}$; brown dash-dot line shows the region obtained via CCP \cite{Shen2022Region}; and red dotted polygon represents the LPF-based approximation \cite{Gan2014linear}.}
\label{fig:CCP_linear_exact}
\end{figure}

The proposed inner approximation region $\mathcal{\tilde{U}}$ is compared with three alternative methods: (1) the linearized power flow (LPF) model in \cite{Gan2014linear}, making a balanced three-phase voltage assumption; (2) the convex-concave procedure (CCP) \cite{Shen2022Region}; and (3) a brute-force search for obtaining the accurate dispatchable region $\mathcal{U}$. All comparisons are conducted under the case \emph{Baseline $\&$ Volt. R}, as defined previously in Section \ref{sec:numerical:subsec:factors}. The comparison results are illustrated in Fig. \ref{fig:CCP_linear_exact}.

Quantitatively, we evaluate each method by measuring the enclosed area relative to the ground-truth region determined by brute-force search. The relative areas covered by the LPF model, CCP algorithm, and our proposed inner approximation $\mathcal{\tilde{U}}$ are 54.1$\%$, 59.2$\%$, and 87.4$\%$, respectively.

Regarding computational efficiency, the adoption of LPF reduces the runtime of Algorithm \ref{Algo:cutting} to 13.2 seconds. However, due to its idealized voltage balance assumption, the LPF-generated region is usually a coarse approximation that is neither strictly inner nor outer of the dispatchable region $\mathcal{U}$.
The CCP method is computationally costly, primarily due to numerical instability arising from ill-conditioned bus injection-based SDP computations \cite{Gan2014linear}. Additionally, its sensitivity to initial point selection further compromises numerical stability and practical reliability.
In contrast, the proposed inner approximation $\mathcal{\tilde{U}}$ attains superior accuracy and computational efficiency. 
Unlike CCP, which requires multiple iterative steps, our approach verifies the feasibility of renewable power $\mathbf{u}$ by solving a single optimization problem. As analytically verified in Corollary \ref{coro:inner:approx} and illustrated in Fig. \ref{fig:CCP_linear_exact}, $\tilde{\mathcal{U}}$ consistently resides within $\mathcal{U}$, ensuring both reliability and practical applicability.

\section{Conclusion}\label{sec:conclusion}
In this paper, we developed both outer and inner approximations of the renewable generation dispatchable region in unbalanced three-phase radial networks. First, we formulated a non-convex optimization problem to define this region, which is then relaxed to a convex SDP.
Using an SDP-based projection algorithm (Algorithm \ref{Algo:cutting}), we accurately derived a polytopic outer approximation of the dispatchable region.
Further, we introduced sufficient conditions to guarantee the exact SDP relaxation by incorporating the power loss as a penalty term, thereby establishing an inner approximation of the dispatchable region.
Our method is adaptable to unbalanced three-phase AC power flow models, providing more precise outcomes compared to prior methods such as linearization of power flow.

In the future work, we aim to develop a polytopic approximation algorithm with faster
computation and a more accurate inner approximation of the dispatchable region.

\input{main.bbl}

\end{document}